\documentstyle[11pt]{article}
\newcommand{\sa}{$\Sigma $ }
\newcommand{\ul}{\underline}

\begin{document}

\title{  $SL(2,C)$-TOPOLOGICAL QUANTUM FIELD THEORY WITH CORNERS}
\author{ R{\u{a}}zvan Gelca}
\maketitle

\begin{abstract} We describe the construction of a topological quantum field
theory with corners for 3-manifolds,
 using  quantum deformations of $sl(2,{\bf C})$.
In the construction there appears 
a sign obstruction for some of the Moore-Seiberg
equations. We solve this problem by means of the Klein four group.
\end{abstract}

\begin{center}
{\bf 0. INTRODUCTION}
\end{center}

\medskip

This paper contains the basic ideas for the construction of an
$sl(2,{\bf C})$ topological quantum field theory for orientable 3-manifolds.
 The construction is done following ideas originated in [RT], where the
authors describe the construction, for any modular Hopf algebra, of
a topological quantum field theory that satisfies the Atiyah-Segal
axioms (a so called smooth topological quantum field theory), see
also [T] and [KM]. 

In [Wa] it is described a construction of a topological quantum field theory 
from a modular Hopf algebra, when one  allows gluings along
surfaces with {\em boundary} (a topological quantum field theory with 
corners). Compared to the construction described in [RT], this approach 
enables us to localize computations, unlike the case of
surgery diagrams
where computations  are global.
This also leads to various ways of computing the Jones polynomial ([J])
for knots.

The key elements in this construction are the decompositions of surfaces into 
disks, annuli and pairs of pants, and the moves that transform
one decomposition into another. A surface together with such a
decomposition is the analogue of  a vector space endowed with
a basis. For topological quantum field theories with corners
this point of view is essential, since in this case the gluings will occur 
along whole subsurfaces on the boundary of 3-manifolds. This is also
a very useful approach even for smooth topological quantum
field theories if one wants to do computations, since the rather 
complicated mapping class group is replaced with the simpler groupoid
of transformations between decompositions.

Unfortunately, $sl(2,{\bf C})$ topological quantum field theory does not fit
in this framework, since there appears a sign problem. This sign
 problem comes from the
fact that whenever one moves a 
coupon of the Weyl element (denoted by $D$ in [KM] and by $w_i$ in
[T]) over a maximum or a minimum 
on a strand labeled by an even dimensional representation, the sign 
of the morphism assigned to the diagram changes. In doing computations
with diagrams, this  produces an obstruction
for some  of the  Moore-Seiberg equations. The problem is also
related to an asymmetry of the skein relation of the Jones polynomial
which will be explained in
[G]. The case where one restricts 
to working only with odd-dimensional irreducible representations has 
been done in [FK].

The sign problem shows that the decompositions into disks, annuli and
pairs of pants are not fine enough to make the modular functor be well
defined. 
In the present paper we solve the sign problem simply by adding
some extra structure, namely by 
introducing a family of simple closed curves indexed by the elements
of the Klein group on the boundary of the 3-manifold. This approach
is similar to that of adding a Lagrangian space ([Wa], [T]) in
order to solve the projective ambiguity of the invariants,
and fits in the formalism of Chap. III in [T]. We deduce the
Moore-Seiberg equations corresponding to our situation and adapt
the techniques from
[FK] to obtain an $sl(2,{\bf C})$-topological
quantum field theory. 
The invariants we  we obtain for closed manifolds are the 
Reshetikhin-Turaev invariants multiplied by $X^{-1}$, where the constant
$X$ is defined in Section 4.

In the case where one wants to avoid the use of Weyl elements, N. 
Reshetikhin suggested us an approach by using a ${\bf Z}_2$ structure on the
boundary of the 3-manifold, similar to our Klein group structure.
This is obtained by orienting (or rather coorienting) the circles that 
decompose the surface into elementary surfaces, which is a 
natural approach if one notes that the oriented graph that is the
core of the surface is the diagram of the associated vector space.
In this way the Moore-Seiberg relations that don't work are replaced
by their double covers, which clearly hold.

The introduction of Weyl elements  makes computations easier and
 gluings more flexible. It also has the advantage that it attaches the
same vector space to homeomorphic surfaces, regardless of the structure,
and the same is true about the morphisms. 

Regarding the formalism of diagrams, let us mention the following.
First, at the level of pairs of pants, there are two nontrivial
relations that have to be satisfied by the modular functor.
 The first one is the
third Reidemeister move. The second one is
the one that shows that the cube of the rotation is the identity,
 and its proof is usually omitted in the
literature. We give a short proof of it. The tetrahedron
of the 6j-symbols has been altered by means of rotations an a sign
change for aesthetic reasons. Also, compared to [FK], our diagrams are
defined in such a way that, when doing computations, all tensor
contractions are performed at the top, while all twistings and
braidings occur at the bottom. This is very well illustrated in
the proof we give of the pentagon identity.

The author wants to thank Ch. Frohman for the fruitful discussions.

\medskip

\begin{center}
{\bf 1. COMPLETELY EXTENDED SURFACES}
\end{center}

\medskip

In this section we introduce the Klein group structure that will make the
modular functor be well defined.
Let $\Sigma $ be a compact, orientable
piecewise linear surface (with boundary).

\medskip

\underline{Definition:} A DAP-decomposition of \sa
consists of 

\ \ - a collection $\alpha $ of disjoint circles in
the interior of \sa that cut \sa into a collection of 
disks, annuli and pairs of pants.

\ \ - an ordering of the components of \sa cut along
$\alpha $.

\ \ - if $\Sigma _0$ is one of these components
then the  boundary components of $\Sigma _0$
should be numbered by 1, 2 and 3 if $\Sigma _0$ is
a pair of pants, 1 and 2 if $\Sigma _0$ is an annulus
and 1 if $\Sigma _0$ is a disk. Each boundary component $C$
of $\Sigma _0$ should have a fixed parametrization 
$f:S^1\rightarrow C$. Fix three disjoint embedded
arcs joining $e^{i\epsilon}$ on the j-th boundary   
component to $e^{-i\epsilon}$ on the (j+1)-st 
boundary component (taken modulo
the number of components of $\Sigma _0$), where $0<\epsilon
<\pi$ is fixed. These arcs will be called seams.

\ \ - if $\Sigma _0$ and $\Sigma _1$ are two of the components of
\sa cut along $\alpha $ and $C$ is a circle in $\Sigma _0\cap \Sigma _1$
then the parametrization of $C$ in $\Sigma _0$ should coincide with the 
complex conjugate of the parametrization of $C$ in $\Sigma _1$.

An example of a DAP-decomposition is given in Fig.1.
Two DAP-decompositions that coincide up to isotopy  will be
considered identical.

\input epsf

\begin{figure}[htbp]
\centering
\leavevmode
\epsfxsize=4.5in
\epsfysize=1.5in
\epsfbox{Quantum/fig1.eps}
 
Fig. 1.1.
\end{figure}

Let $K=\{e,a,b,c\}$ be the Klein four group. Let us recall
that the elements of $K$ satisfy $a^2=b^2=c^2=e$, $ab=c$,
$bc=a$, $ac=b$.

\medskip

\underline{Definition:} A completely extended surface
is a triple $(\Sigma ,D,B)$ where

\ - \sa is a compact, oriented, piecewise linear surface,

\ - $D$ is a DAP-decomposition,

\ - $B$ is a collection of bands (circles) such that for each
elementary surface $\Sigma _0$ in the  DAP-decomposition of
\sa there is one band parallel to each of the boundary
components of $\Sigma _0$. The bands are indexed by $e, a,b$,
or $c$ according to the following rules:

i) if the two numbers that correspond to the 
numbering of a decomposition circle in two neighboring elementary 
surfaces are both 1, or none of them is 1, then the product
of the indices of the two bands parallel to that component should
have the product equal to $a$ or $b$;

 ii) in the other cases the product of those indices
should be $e$ or $c$.

\medskip

The reason why we have introduced the bands instead of
simply indexing the decomposition circles is that we want to make gluings
with corners more flexible. 
An example, with the seams omitted is given in Fig 1.2.

\begin{figure}[htbp]
\centering
\leavevmode
\epsfxsize=4.4in
\epsfysize=1.5in
\epsfbox{Quantum/fig2.eps}

Fig. 1.2.
\end{figure}

In the picture above the seams have been omitted. We will do this
whenever it will seam unlikely to cause confusion.

For simplicity we will make the notations described in Fig 1.3.

\begin{figure}[htbp]
\centering
\leavevmode
\epsfxsize=4.2in
\epsfysize=1.4in
\epsfbox{Quantum/fig3.eps}

Fig. 1.3.
\end{figure}

We will factor the set of completely extended surfaces (shortly
ce-surfaces) by the following two identifications:

1. Is described in Fig 1.4, where $\bar{u}=u$ if $u=e$ or $c$ and
$\bar{u}=cu$ if $u=a$ or $b$.

\begin{figure}[htbp]
\centering
\leavevmode
\epsfxsize=3.1in
\epsfysize=0.8in
\epsfbox{Quantum/fig4.eps}

Fig. 1.4.
\end{figure}

2. If $(\Sigma ',D',B')\subset (\Sigma ,D, B)$ is a subsurface,
let us consider the ce-surface obtained from $(\Sigma ,D, B)$
by multiplying all bands
along the boundary components of $\Sigma '$ by $c$. Identify the
newly obtained ce-surface with $(\Sigma ,D, B)$. 

These identifications can be easily understood in the ``arrow''
notation. The first one means that we are allowed to slide bands
over a decomposition circle, like in Fig. 1.5, while the second means that
we are allowed to reverse one arrow on each band along the boundary of a 
subsurface of \sa.

\begin{figure}[htbp]
\centering
\leavevmode
\epsfxsize=3.2in
\epsfysize=0.8in
\epsfbox{Quantum/fig5.eps}

Fig. 1.5.
\end{figure}

 For example, the surface in Fig. 1.6 is obtained from Fig. 1.2 by using the
first kind of identification. This is further equal to the one from Fig. 1.7
by performing the second kind of identification on the dotted surface.

\begin{figure}[htbp]
\centering
\leavevmode
\epsfxsize=4.5in
\epsfysize=1.5in
\epsfbox{Quantum/fig6.eps}

Fig. 1.6.
\end{figure}

\begin{figure}[htbp]
\centering
\leavevmode
\epsfxsize=4.5in
\epsfysize=1.5in
\epsfbox{Quantum/fig7.eps}

Fig 1.7.
\end{figure}

\underline{Remark} We do not allow
bands to slide from one end of an annulus to the other.

\medskip

At this point the construction looks very artificial. It will
become very natural when we define the functor, the band structure will
mimic the behavior of the coupons labeled by $D$, the arrows above
will reflect the direction of the arrows on a strand over a maximum or
a minimum.

From now on the structure on a surface that makes it into a ce-surface
will be called a DB-structure. We have come to the point where we
introduce a set of elementary moves that transform one DB-structure
into another.
At the level of the modular functor, they will correspond to
changes of basis in the associated vector space.
These moves are $K$, which  is the multiplication by $c$ on one band
(see Fig 1.8), $T_1$, which is the twist given in Fig. 1.9, where the
band structure is left invariant, $B_{23}$ described in Fig. 1.10,
leaving also the band structure invariant, the rotation $R$ 
consisting of the permutation
of the numbers 1, 2 and 3, multiplication by $a$ on the top boundary 
component, and by $b$ on the bottom right boundary component (see Fig 1.11.)
 The moves $F$, $S$, $A$ and $D$ are described in the figures next four
figures. 

\begin{figure}[htbp]
\centering
\leavevmode
\epsfxsize=3.2in
\epsfysize=0.8in
\epsfbox{Quantum/fig8.eps}

Fig. 1.8.
\end{figure}

\begin{figure}[htbp]
\centering
\leavevmode
\epsfxsize=3.2in
\epsfysize=1.1in
\epsfbox{Quantum/fig9.eps}

Fig. 1.9.
\end{figure}

\begin{figure}[htbp]
\centering
\leavevmode
\epsfxsize=3.2in
\epsfysize=1.1in
\epsfbox{Quantum/fig10.eps}

Fig. 1.10.
\end{figure}

\begin{figure}[htbp]
\centering
\leavevmode
\epsfxsize=3.2in
\epsfysize=1.1in
\epsfbox{Quantum/fig11.eps}

Fig. 1.11.
\end{figure}

\begin{figure}[htbp]
\centering
\leavevmode
\epsfxsize=3.2in
\epsfysize=1.1in
\epsfbox{Quantum/fig12.eps}

Fig. 1.12.
\end{figure}

\begin{figure}[htbp]
\centering
\leavevmode
\epsfxsize=3.2in
\epsfysize=1.1in
\epsfbox{Quantum/fig13.eps}

 Fig. 1.13.
\end{figure}

\begin{figure}[htbp]
\centering
\leavevmode
\epsfxsize=3.1in
\epsfysize=1.5in
\epsfbox{Quantum/fig14.eps}

Fig. 1.14.
\end{figure}

\begin{figure}[htbp]
\centering
\leavevmode
\epsfxsize=3.1in
\epsfysize=1.1in
\epsfbox{Quantum/fig16.eps}

Fig. 1.15.
\end{figure}

A composition of the elementary moves described above will be called
a move.

\underline{Notes} 1. The moves $F$ and $S$ can only
be applied when the arrows point in the prescribed direction.
As a matter of fact we have no elementary move for the 
$S$-move on the torus, however this can be obtained 
as a composition of moves as seen in Fig. 1.16.

\begin{figure}[htbp]
\centering
\leavevmode
\epsfxsize=4.6in
\epsfysize=2.7in
\epsfbox{Quantum/fig17.eps}

Fig. 1.16.
\end{figure}

2. One should not confuse the moves $T_1, B_{23}, R$ and
$S$ with the corresponding maps between surfaces. Here
they only change the DB-structure, their underlying 
homeomorphism is the identity.

\begin{figure}[htbp]
\centering
\leavevmode
\epsfxsize=3.7in
\epsfysize=1.1in
\epsfbox{Quantum/fig18.eps}

Fig. 1.17.
\end{figure}

3. We factor everything by the move described in Fig 1.17,
 so we make the 
identification from Fig. 1.18, giving the rule of sliding a band in an annulus.

\begin{figure}[htbp]
\centering
\leavevmode
\epsfxsize=1.8in
\epsfysize=1in
\epsfbox{Quantum/fig19.eps}

Fig. 1.18.
\end{figure}

Let us remark that the image of a DB-structure through any of the
elementary moves described above is again a DB-structure, satisfying
the conditions from the definition. 

Let $C$ be a boundary component of \sa. We say that $C$ is
of positive type $(+)$ if it is numbered by a 1 and its associated
band is indexed by an $e$ or a $c$, or if it is numbered by 2 or 3
and its associated band is indexed by an $a$ or a $b$.
In the other cases we say that $C$ is of negative type $(-)$.

\medskip

\underline{Definition:} (Gluing ce-surfaces) Let $(\Sigma ,D,B)$
be a completely extended surface. Let $C_1$ and
$C_2$ be two boundary components of \sa of opposite type
and $f:C_1\rightarrow C_2$ be the homeomorphism induced by the
complex conjugation on $S^1$. Define $(\Sigma _f, D_f, B_f)$ 
to be the ce-surface with $\Sigma _f$ the gluing of \sa by $f$,
$D_f$ and $B_f$ being the DAP-decompositions, respectively band
structures induced by $D$ and $B$.

For example the surfaces from Fig. 1.19 can be glued together but
the ones from Fig. 1.20 cannot. 

\begin{figure}[htbp]
\centering
\leavevmode
\epsfxsize=3in
\epsfysize=0.8in
\epsfbox{Quantum/fig20.eps}

Fig. 1.19.
\end{figure}

\begin{figure}[htbp]
\centering
\leavevmode
\epsfxsize=2in
\epsfysize=0.7in
\epsfbox{Quantum/fig21.eps}

Fig. 1.20.
\end{figure}

\underline{Remark} The condition that $C_1$ and $C_2$ are
of opposite types guarantees that the newly obtained surface 
satisfies the conditions from the definition.

\medskip

\underline{Definition:} A ce-morphism is a map between two ce-surfaces
\begin{eqnarray*}
(f,n):(\Sigma _1,D_1,B_1)\rightarrow (\Sigma _2, D_2, B_2)
\end{eqnarray*}
where $f$ is a homeomorphism satisfying the property that for any 
boundary component $C$ of $\Sigma _1$, $C$ and $f(C)$ are of the same
type, and $n$ is an integer.

\medskip

\underline{Remarks} 1. The 
moves $K$, $T_1$, $B_{23}$, $R$, $F$, $S$, $A$ and $D$
will be considered as ce-morphisms with the underlying homeomorphism
equal to the identity, and paired  with the integer $n=0$.

2. The definition above shows that there
is no morphism between 
the two surfaces in Fig. 1.21.

\begin{figure}[htbp]
\centering
\leavevmode
\epsfxsize=3.5in
\epsfysize=1.1in
\epsfbox{Quantum/fig22.eps}

Fig. 1.21.
\end{figure}

\noindent It follows that the mapping class groupoid of the ce-surfaces having
the same underlying surface is not always connected, in fact it has 
exactly as many components as the number of the boundary components of 
the surface. However, this produces no restrictions on the
gluings along surfaces with boundary, since by an operation of type
$A$ we can always adjust the surface such that the gluing can be done.
The groupoid is connected if the surface has no boundary.

\medskip

\underline{Proposition 1.1.} If $(f,n):(\Sigma _1,D_1,B_1)
\rightarrow (\Sigma _2, D_2, B_2)$
is a ce-morphism then there exists a move between $(\Sigma _2,f(D_1),f(B_1))$
and $(\Sigma _2, D_2, B_2)$. In particular $(id,0):(\Sigma ,D_1,B_1)
\rightarrow(\Sigma ,D_2,B_2)$ is just a move.

\medskip

{\underline{Proof:}}  Since any two DAP-decompositions can be
transformed one into the other by a move, we may assume that $f(D_1)=D_2$.
From the definition of a ce-surface, the indices of the bands of $f(B_1)$
that lie in the interior of $\Sigma _2$ agree with those of $B_2$
modulo multiplication by $c$, so they can be made to agree by
performing some moves of type $K$. From the definition of ce-morphisms
it follows that the same is true for the bands along the boundary
components, and the conclusion follows.$\Box$ 

%\vspace{4in}
\begin{figure}[htbp]
\centering
\leavevmode
\epsfxsize=2.9in
\epsfysize=3.2in
\epsfbox{Quantum/fig23.eps}

Fig. 1.22.
\end{figure}

\underline{Definition:} (Composition of ce-morphisms) If 
$(f_1,n_1):(\Sigma _1,D_1,B_1)
\rightarrow (\Sigma _2, D_2, B_2)$ and $(f_2,n_2):(\Sigma _2,D_2,B_2)
\rightarrow (\Sigma _3, D_3, B_3)$ are two ce-morphisms we define
\begin{eqnarray*}
(f_2,n_2)(f_1,n_1):=(f_2f_1,n_2+n_1-\sigma((f_2f_1)_*L_1,(f_2)_*L_2,L_3))
\end{eqnarray*}
where $L_i$ is the subspace of $H_1(\Sigma _i)$ generated by 
$D_i$ (including the curves on the boundary), $i=1,2,3$
and $\sigma $ is Wall's non-additivity function ([W]).

The negative sign in front of the $\sigma$ function appears because
the cocycle that gives the extension comes from the signature of a 4-manifold,
hence is the negative of Wall's $\sigma$ function.

\underline{Definition:} The dual of a ce-surface $(\Sigma ,D,B)$
is defined to be $-(\Sigma ,D,B):=(-\Sigma ,-D,-B)$ where $-\Sigma :=$\sa
with opposite orientation, $-D$ is obtained from $D$ by
 applying the orientation
reversing moves described in Fig. 1.22 to each of the disks, annuli
or pairs of pants from the decomposition of \sa
and $-B$ is defined as follows:

\ \ i) inside \sa the band structure remains unchanged;

\ \ ii) we multiply by $a$ the index of all bands along boundary components
of positive type, and by $b$ the index for boundary components of
negative type.

\medskip

An example is given in 
Fig. 1.23.

\begin{figure}[htbp]
\centering
\leavevmode
\epsfxsize=4.3in
\epsfysize=1.5in
\epsfbox{Quantum/fig24.eps}

Fig. 1.23.
\end{figure}

\medskip

\underline{Proposition 1.2.} The operation of taking the dual is natural
with respect to morphisms and gluings, and is reflexive.

\medskip

\ul{Proof:} Naturality is straightforward. Reflexivity follows from
the fact that we identify a ce-surface with the ce-surface obtained by
multiplying the indices of all bands along boundary components
by $c$.

\bigskip

\begin{center}
{\bf 2. COMPLETELY EXTENDED 3-MANIFOLDS}
\end{center}
\bigskip

In this section we introduce the manifolds whose boundaries are ce-surfaces.

\underline{Definition:} We say that $(M, D, B, n)$ is a completely extended
3-manifold (shortly ce-3-manifold) if $M$ is a compact, oriented, piecewise 
linear 3-manifold, $D$ is a DAP-decomposition of the boundary of
$M$, $B$ is a collection of bands indexed by the elements
of the Klein group on the boundary of $M$
such that $D$ and $B$ determine a DB-structure on the boundary,
and $n$ is an integer (usually called framing).

\medskip

\underline{Definition:} (of the boundary) $\partial (M,D,B,n)=
(\partial M,D,B)$.

\medskip

\ul{Definition:} (disjoint union) $(M_1,D_1,B_1,n_1)\sqcup (M_2,D_2,B_2,n_2)=
(M_1\sqcup M_2,D_1\sqcup D_2, B_1\sqcup B_2, n_1+n_2)$.

\medskip

\underline{Definition:} (the mapping cylinder) If $(f,n):(\Sigma _1, D_1,
B_1)\rightarrow(\Sigma _2, D_2,B_2)$ is a ce-morphism, we define the
mapping cylinder of $(f,n)$ to be $I_{(f,n)}:=(M, D, B, n)$
where $M=I_f$ the mapping cylinder of $f$ and $(\partial M, D,B)=
-(\Sigma _1,D_1,B_1)\cup (\Sigma _2, D_2,B_2)$, where the two surfaces
are glued together along the boundary components that are mapped one into
the other by $f$.

\medskip

Here we consider the definition of the mapping cylinder in which the 
boundary of the first surface is identified with the boundary of the second.

\medskip

\underline{Definition:} (gluing of ce-3-manifolds) Let $(M,D,B,n)$ be a 
ce-3-manifold. Let $(\Sigma _i, D_i,B_i)\subset \partial (M,D,B,n)$, $i=1,2$
be two disjoint ce-surfaces.
Let $(f,m):(\Sigma _1, D_1,
B_1)\rightarrow -(\Sigma _2, D_2,B_2)$ be a ce-morphism. Define the gluing
of $(M,D,B,n)$ by $(f,m)$ to be 
\begin{eqnarray*}
(M,D,B,n)_{(f,m)}:=(M_f,D',B',m+n-\sigma (K,L_1\oplus L_2, \Delta ^-))
\end{eqnarray*}
where

- $M_f$ is the (piecewise linear) gluing of $M$ by $f$,

- $D'$  is the image of $D$  under the quotient map
$(\partial M\backslash Int(\Sigma _1\cup \Sigma _2))\rightarrow \partial M_f$,

- $B'$ is obtained from $B$ by only keeping the bands that are contained in
the complement of the two surfaces that we glued along,

- $\sigma $ is Wall's non-additivity function (see [W]),

- $K$ is the kernel of $H_1(\Sigma _1\cup \Sigma _2)
\rightarrow H_1(M)/J$, where $J$ is the subspace of $H_1(\partial M)$
spanned by the decomposition curves lying in the complement of $int(\Sigma
_1\cup \Sigma _2)$,

- $L_i$ are the subspaces of $H_1(\Sigma _i)$ generated by $D_i$, $i=1,2$,

- $\Delta ^-:=\{(x,-f_*(x))$ $|x\in H_1(\Sigma _1)\}$.

\medskip

For a better understanding let us briefly review the explanation from [Wa]
concerning this choice of the framing.

Let $L$ be the subspace of $H_1(\partial M)$ generated by $D$. 
Let $M^+$ be a 3-manifold such that $\partial M^+=-\Sigma$
and $Ker(H_1(\partial M)\rightarrow H_1(M^+))=L$. Let $W$ be e 4-manifold 
bounded by $M\cup M^+$ (glued along $\partial M$) whose signature is $n$.

Now let us cap $\Sigma _1$ and $\Sigma _2$ with disks along the boundary
components, disks lying entirely in $M^+$, and obtain two surfaces 
$\hat{\Sigma _1}$ and $\hat{\Sigma _2}$ in $M^+$. One can assume that these two
surfaces bound two disjoint manifolds in $M^+$, say $N_1$ and $N_2$
(see Fig. 2.1.)

\begin{figure}[htbp]
\centering
\leavevmode
\epsfxsize=3in
\epsfysize=1.2in
\epsfbox{Quantum/fig25.eps}

Fig. 2.1.
\end{figure}

If we consider $W'$ the 4-manifold bounded by $N_1\cup I_f\cup N_2$,
with signature $m$, where $I_f$ is the mapping cylinder of $f$
extended to the disks that we capped with, then we see that the framing
of $(M,D,B,n)_{(f,m)}$ is the signature of the 4-manifold obtained by
gluing $W$ and $W'$.

\medskip

\underline{Remark} Since the morphism $(f,m)$ is defined from
one surface to the dual of the other, we see that the boundary of the newly
obtained ce-3-manifold satisfies the conditions from the definition of
a ce-surface.

\bigskip

\underline{Proposition 2.1.} The following properties hold:

1) the gluing operation is associative,

2) $I_{(g,n)(f,m)}=I_{(g,n)}\cup _{(id,0)}I_{(f,m)}$,

3) $(M,D,B,n)_{(f,m)}=(M,D,B,n)\cup _{(id,0)}I_{(f,m)}$.

\medskip

\ul{Proof:} 1) The only thing one has to prove is that the value of the
framing does not depend on the order in which we perform the gluings.
However, this follows by interpreting the framing as the signature
of a 4-manifold as it has been done above, since the gluing of
manifolds is associative. 

2) and 3) also follow from a similar argument.$\Box$

\bigskip

\begin{center}
{\bf 3. THE AXIOMS OF A TOPOLOGICAL QUANTUM FIELD THEORY FOR CE-SURFACES 
AND 3-MANIFOLDS}
\end{center}
\bigskip

In this section we are going to describe the axioms that a topological quantum
field theory with corners for ce-surfaces and 3-manifolds should satisfy.
 They are identical with  those introduced by Walker [Wa]
for extended manifolds.

A label set is a finite set ${\cal L}$ equipped with an involution
$x\rightarrow \hat{x}$ and a 
distinguished element $1\in {\cal L}$, with $\hat{1}=1$.
Define the category of labeled completely extended surfaces (lce-surfaces),
whose objects are ce-surfaces with a label attached to each boundary
component, and whose morphisms are ce-morphisms which preserve the
labeling. From now on we will denote ce-surfaces and ce-3-manifolds
by a single letter, usually $\sigma $ for surfaces and $\mu $ for
3-manifolds.
Thus a lce-surface is a pair $(\sigma , l)$ where $\sigma$ is
a ce-surface and $l$ is a function from the set of boundary
components of $\sigma$ to ${\cal L}$.

\medskip

\underline{Definition:} A topological quantum field theory of
label set ${\cal L}$ consists of 

\ i) a functor $V$ from the category of lce-surfaces to the
category of finite dimensional vector spaces and morphisms;

\ ii) an assignment $\mu \rightarrow Z(\mu )\in V(\partial \mu)$ for each
ce-3-manifold $\mu $;

satisfying the axioms below:

\medskip

(3.1.) Disjoint union axiom

\ \ $V(\sigma _1 \sqcup \sigma _2,l_1\sqcup l_2)=V(\sigma _1,l_1)\otimes 
V(\sigma _2, l_2)$.

\medskip

(3.2.) Gluing axiom for $V$

Let $\sigma $ be a ce-surface, $C$ and $C'$ two sets of boundary
components of $\sigma $, $g:C\rightarrow C'$ the parametrization reflecting
map. Let $\sigma _g$ be $\sigma $ glued by $g$. Then
\begin{eqnarray*}
V(\sigma _g,l)=\bigoplus _{x\in {\cal L}(C)}V(\sigma ,(l,x,\hat{x}))
\end{eqnarray*}
where the sum is over all labelings of $C$.

\medskip

(3.3.) Duality axiom
\begin{eqnarray*}
V(\sigma ,l)=V(-\sigma,\hat{l})^*
\end{eqnarray*}
satisfying the following compatibility conditions 

\medskip

-the identifications 

\centerline{$V(\sigma ,l)=V(-\sigma,\hat{l})^*$}

\centerline{$V(-\sigma,\hat{l})=V(\sigma ,l)^*$}

are mutually adjoint;

\medskip

-if $\phi =(f,n):(\sigma _1,l_1)\rightarrow (\sigma _2,l_2)$,
let $\stackrel{-}{\phi}:=(f,-n):(-\sigma _1,\hat{l}_1)
\rightarrow (-\sigma _2,\hat{l}_2)$.
Then for any $x\in V(\sigma _1,l_1)$ and $y\in V(-\sigma _1,\hat{l}_1)$
we have
\begin{eqnarray*}
<x,y>=<V(\phi )x, V(\stackrel{-}{\phi})y>
\end{eqnarray*}
where $<,>$ is the pairing between the space and its dual.
This condition says that $V(\phi )$ is the adjoint inverse of 
$V(\stackrel{-}{\phi})$;

\medskip

-if $\alpha _1\otimes \alpha _2\in V(\sigma _1\sqcup \sigma _2)=
V(\sigma _1)\otimes V(\sigma _2)$

\ $\beta _1\otimes \beta _2\in V(-\sigma _1\sqcup -\sigma _2)=
V(-\sigma _1)\otimes V(-\sigma _2)$

\noindent then 

$<\alpha _1\otimes \alpha _2,\beta _1\otimes \beta _2>=<\alpha _1,\beta_1>
<\alpha _2, \beta _2>$;

\medskip

-there is a function $S:{\cal L}\rightarrow {\bf C}^*$ such that
if 

\ $\oplus _{x\in {\cal L}(C)}\alpha _x\in V(\sigma _g, l)=\bigoplus_{
x\in {\cal L}(C)}V(\sigma ,(l,x, \hat{x}))$

\ $\oplus _{x\in {\cal L}(C)}\beta _x\in V(-\sigma _g, \hat{l})=\bigoplus_{
x\in {\cal L}(C)}V(-\sigma ,(\hat{l}, \hat{x},x))$

\noindent then

\ $<\oplus _x \alpha _x,\oplus _x \beta _x>=\Sigma _x S(x)<\alpha _x, 
\beta _x>$

\noindent where $C\subset  \partial \sigma $ consists out of
$n$ circles, $x=(x_1,x_2,\cdots ,x_n)$ and $S(x)=S(x_1)S(x_2)
\cdots S(x_n)$.

\medskip

(3.4.) Empty surface axiom

\ \ $V(\emptyset )\tilde{=}{\bf C}$.

\medskip

(3.5.) Disk axiom

 If $D$ is a ce-disk we have $V(D,m)={\bf C}$ if $ m=1$ and $0$ otherwise.

\medskip

(3.6.) Annulus axiom

If $A$ is a ce-annulus we have $V(A,(m,n))={\bf C}$ if $ m=\hat{n}$ and $0$ 
otherwise.

\medskip

(3.7.) Disjoint union axiom for $Z$

\centerline{$Z(\mu _1 \sqcup \mu _2)=Z(\mu _1)\otimes Z(\mu _2).$}

\medskip

(3.8.) Naturality axiom

Let $(f,m):(M_1,D_1,B_1,n_1)\rightarrow (M_2,D_2,B_2,n_2)$
be a morphism, and suppose that $n_2=n_1+m+\sigma (K,L_1,L_2)$,
where $K=ker(H_1(\partial M_2)\rightarrow H_1(M_2))$, and $L_1$ and
$L_2$ are the (Lagrangian) subspaces of $H_1(\partial M_2)$
generated by $f(D_1)$ and $D_2$.
Then 
\begin{eqnarray*}
V(f|\partial (M_1,D_1,B_1,n),m)Z(M_1,D_1,B_1,n_1)=
Z(M_2,D_2,B_2,n_2).
\end{eqnarray*}

\medskip

(3.9.) Gluing axiom for $Z$

Let $\mu$ be a ce-3-manifold, $\sigma _1, \sigma _2\subset \partial \mu$
two disjoint ce-surfaces, and suppose that there exists a ce-morphism
$\phi :\sigma _1\rightarrow -\sigma _2$. Then one can define $\mu _{\phi}$,
the gluing of $\mu $ by $\phi $. We have
\begin{eqnarray*}
V(\partial \mu )=\bigoplus _{l_1,l_2} V(\sigma _1,l_1)\bigotimes
V(\sigma _2,l_2)\bigotimes V(\partial \mu \backslash (\sigma _1 \cup
\sigma _2),(\hat{l}_1,\hat{l}_2))
\end{eqnarray*}
where $l_i$ runs through all labelings of $\sigma _i$. Also
\begin{eqnarray*}
Z(\mu )=\bigoplus _{l_1,l_2}\Sigma _j\alpha _{l_1}^j\otimes \beta _{l_2}^j
\otimes \gamma _{\hat{l}_1\hat{l}_2}^j.
\end{eqnarray*}

The axiom states that 
\begin{eqnarray*}
Z(\mu _{\phi})=\bigoplus _l \sum_j <V(\phi)\alpha_l^j,\beta_{\hat{l}}^j>
\gamma _{l\hat{l}}^j
\end{eqnarray*}
where $l$ runs through all labelings of $\sigma _1$.

\medskip

(3.10.) Mapping cylinder axiom

For $(id,0):\sigma \rightarrow \sigma $ we have
\begin{eqnarray*}
V(\partial I_{(id,0)})=\bigoplus _l V(\sigma ,l)\otimes V(\sigma ,l)^*.
\end{eqnarray*}

If $id_l$ is the identity matrix in $V(\sigma ,l)\otimes V(\sigma ,l)^*$
then 
\begin{eqnarray*}
Z(I_{(id,0)})=\bigoplus _{l\in {\cal L}(\sigma )}id_l
\end{eqnarray*}

\bigskip

\underline{Remark}  As a consequence of 
the gluing axiom, the mapping cylinder axiom has the 
following stronger form. Let $\phi :\sigma _1\rightarrow \sigma _2$
be a ce-morphism. Then $Z(I_{\phi})=\bigoplus _{l\in {\cal L}(\sigma _1)}
V(\phi _l)$ where $\phi _l:(\sigma _1,l)\rightarrow (\sigma _2,l)$
are the corresponding lce-morphisms.

\bigskip
\begin{center}
{\bf 4. THE CONSTRUCTION OF THE BASIC DATA}
\end{center}
\bigskip

A basic data for a topological quantum field theory
with label set ${\cal L}$ consists of:

\  \ - an assignation of vector spaces to labeled, completely
extended disks, annuli and pairs of pants;

\ \ - a choice of basis elements in these vector spaces;

\ \ - a definition of the morphisms of the vector spaces of lce-disks,
annuli and pairs of pants induced by the operation of taking the dual,
which identifies each vector space with its dual. This identification
should come together with a pairing between the space and its dual;

\ \ - a description of the operators that correspond to the moves
$K_i$, $T_i$, $R$, $B_{23}$, $F$, $S$, $P$, $A$, $D$, $C=(id,1)$ and $\psi$;

\ \ - a function $S:{\cal L}\rightarrow {\bf C}^*$ (needed in axiom (3.3.)).

\medskip

Using the 10 axioms from  the previous section we see that the basic data
uniquely defines $V$ and $Z$, provided that some compatibility conditions
are satisfied.

We will construct a basic data by means of quantum deformations of
the universal enveloping algebra of $sl(2,{\bf C}) $. Let us remind some
 basic facts about the algebra ${\cal A}_r$ introduced by Reshetikhin and
Turaev in [RT] (see also [KM]). Note that this algebra is denoted by
$U_t$ in [RT], our notation is the one from [KM].

Let $r$ be an integer,  $s=e^{\pi i /r}$ and 
$t=e^{\pi i/2r}$. For an integer $n$ one defines 
\begin{eqnarray*}
[n]=\frac{s^n-\bar{s}^n}{s-\bar{s}}=\frac{sin\frac{\pi n}{r}}
{sin\frac{\pi}{r}}
\end{eqnarray*}
$[n]!=[1][2]\cdots [n], \ \ [0]!=1$,
\begin{eqnarray*}
\left[
\begin{array}{clcr}
n \\
k
\end{array}
\right] 
=\frac{[n]!}{[k]![n-k]!},
\ X=\sqrt{\sum_{k=1}^{r-1}[k]^2}=\frac{\sqrt{\frac{r}{2}}}{sin\frac{\pi}{r}}
\end{eqnarray*}
The number $X$ is denoted by ${\cal D}$ in [T] and is $1/b$ in
[KM] 
(one should not confuse this $X$ with the operator $X$ defined below).

${\cal A}_r$ is the associative algebra over ${\bf C}$, with unit, generated
by $X,Y, K$ and $\bar{K}$ satisfying the following relations:

$K\bar{K} =\bar{K}K=1$, $K^{4r}=1$, $X^r=Y^r=0$, $KX=sXK$,
$KY=\bar{s}YK$, $XY-YX=(K^2-\bar{K}^2)/(s-\bar{s})$.

\medskip

Let us recall that ${\cal A}_r$  can be given
a Hopf algebra structure (see [RT] or [KM]).

Let $k$ be an integer $1\leq m<r$, and let $m=(k-1)/2$. Let $\ul{k}$ be the
vector space of dimension $k$, spanned by the vectors $\{e_j\}$ where
$j$ runs from $-m$ to $m$ with step 1. Consider the action of ${\cal A}_r$
on $\ul{k}$ defined by

\ \ $Xe_j=[m+j+1]e_{j+1}, \ j<m, \ Xe_m=0$,

\ \ $Ye_j=[m-j+1]e_{j-1}, \ j>-m, \ Ye_{-m}=0$,

\ \ $Ke_j=t^{2j}e_j$.

By [RT], $\ul{k}$ is an irreducible representation of ${\cal A}_r$. 

Let us also 
recall the definitions of 
two operators that we need for our construction. 
The first one is the Weyl element that 
provides an ${\cal A}_r$-isomorphism between $\ul{k}$
and  $\ul{k}^*$, where $\ul{k}^*$is the 
dual representation of $\ul{k}$, and is given by  
\begin{eqnarray*}
D(e^j)=\left[
\begin{array}{clcr}
2m \\
m-j
\end{array}
\right]   
(it)^{2j}e_{-j}
\end{eqnarray*}
\noindent where $\{e^j\}$ is the
basis dual to $\{e_j\}$ 
One should note that $D^*D^{-1}=(-1)^{k-1}K^2$ (see [KM]). The $(-1)^{k-1}$
that appears in this formula is responsible for the obstruction of constructing
the topological quantum field theory, and made the introduction of the bands 
on the surface necessary.

The second operator is ${\mathaccent20{R}}:\ul{k}\otimes \ul{k}'\rightarrow
\ul{k}'\otimes \ul{k}$ given by
\begin{eqnarray*}
\mathaccent20{R}(e_i\otimes e_j)=\sum_{n\geq 0, j+n\leq m,j-n\geq -m'}
\frac{(s-\bar{s})^n}{[n]!}\frac{[m+i+n]!}{[m+i]!}\frac{[m'-j+n]!}{[m'-j]!}
t^{4ij-2n(i-j)-n(n+1)}e_{j-n}\otimes e_{i+n}.
\end{eqnarray*}
where $m=(k-1)/2$ and $m'=(k'-1)/2$

Following [FK] we define a non-degenerate sesquilinear form
$($ , $):\ul{k}\otimes \ul{k}\rightarrow {\bf C}$. If $v\in \ul{k}$,
$v=\Sigma a_je_j$ define $\bar{v}=\Sigma \bar{a}_je_j$.
Then
\begin{eqnarray*}
(w,v)=i^{k-1}D^{-1}(w)(\bar{v}).
\end{eqnarray*}
Extend $($ , $)$ to tensor products by 
\begin{eqnarray*}
(w_1\otimes w_2,v_1\otimes v_2)=(w_1,v_1)(w_2,v_2)
\end{eqnarray*}

For $\alpha $ a map between tensor products of representations,
let us denote by $\alpha ^*$ its adjoint with respect to this
pairing. If $\alpha ^*\alpha =1$ we call $\alpha $ an isometry;
an invertible isometry is called unitary. As an example $\mathaccent20{R}$ is
unitary [FK].

\bigskip

\ul{Definition} [Wa]: A representation $V$ of ${\cal A}_r$ is $bad$
if for any ${\cal A}_r$-linear map $\alpha :V\rightarrow V$ the trace of
$K^2\alpha $ is zero.

\medskip

By Theorem 8.4.3 in [RT], if $\ul{m}$ and $\ul{n}$ are as before
\begin{eqnarray}
\ul{m}\bigotimes \ul{n}=\bigoplus_p\ul{p}\bigoplus B
\end{eqnarray}
where $B$ is bad, and the sum is taken over all $p$ with
$|m-n|+1\leq p \leq min\{m+n-1,2r-2-m-n\}$, and $m+n+p$ odd,
and the decomposition is unique.

As remarked in [Wa], if $V$ is bad and $W$ is any representation then 
$W\bigotimes V$ is also bad. This enables us to restrict ourselves only
to the part of the tensor product $\ul{m}\otimes \ul{n}$ which
is good, i.e. a sum of irreducible representations, and from now on,
whenever we write $\ul{m}\otimes \ul{n}$ we will actually mean the good part.
It is of no difficulty to check that the tensor product defined this way
satisfies all the properties of the usual tensor product. The morphisms
also behave nicely under taking the quotient over the bad part.

A class of spaces that will be of interest in the sequel are the spaces of
${\cal A}_r$-linear maps $\alpha :\ul{m}\otimes \ul{n}\rightarrow \ul{p}$,
which will be denoted by $V^{mn}_p$, where $\ul{m},\ul{n}$ and
$\ul{p}$ are irreducible representations.
By (1) and Schur's lemma, these spaces are either one or
zero dimensional.
As an example if $k$ is an integer and $m=(k-1)/2$ then
$V^{kk}_1$ is generated  by the orthogonal 
projection $\psi _k:\ul{k}\otimes 
\ul{k}\rightarrow \ul{1}$ given by
\begin{eqnarray*}
\psi _(e_j\otimes e_h)=\delta _{j,-h}\frac{(-1)^{2m}}{\sqrt{[2m+1]}}
\left[
\begin{array}{clcr}
2m \\
m-j
\end{array}
\right]
(it)^{-2j}e_0
\end{eqnarray*}
where $\delta $ is the Kronecker symbol.
\medskip
In the future, the  ${\cal A}_r$-linear orthogonal projections
will be simply called projections, their adjoints will be called 
inclusions.
This terminology is very natural if one considers relation (1).

Define the (normalized) trace pairing 
$<$ , $>_t:V^{mn}_p\bigotimes V^{mn}_p\rightarrow
{\bf C}$ by
\begin{eqnarray*}
<\alpha ,\beta >_t\cdot 1_p=\alpha \circ  \beta ^*
\end{eqnarray*}
where $1_p:\ul{p}\rightarrow \ul{p}$ is the identity. Schur's lemma
implies that the trace pairing is well defined. In the case when we consider
only odd-dimensional representations, it has been proved in [FK]
that the trace pairing is positive definite. The proof can be
adapted to show that the same is true in the case when we also
allow even-dimensional representations. It follows that the trace
pairing admits orthogonal projections, namely maps $\alpha \in V_p^{mn}$ 
with $<\alpha ,\alpha >_t=1$.

We assume that the reader is familiar with the formalism of diagrams 
described in  [RT], [KM] and [T], when working with morphisms of ${\cal A}_r$
representations. In addition, the empty coupons will stand for the morphism 
$D$ or its inverse.
For morphisms $\alpha :\ul{m}\otimes \ul{n}\rightarrow
\ul{p}$ (i.e. $\alpha \in V^{mn}_p$)  and $\beta :\ul{p}\rightarrow
\ul{m}\otimes \ul{n}$ we make the notation from Fig. 4.1.

\begin{figure}[htbp]
\centering
\leavevmode
\epsfxsize=2.7in
\epsfysize=0.8in
\epsfbox{Quantum/fig31.eps}

Fig. 4.1.
\end{figure}

 The trace pairing is now described by the diagram in Fig. 4.2.

\begin{figure}[htbp]
\centering
\leavevmode
\epsfxsize=1.8in
\epsfysize=1.3in
\epsfbox{Quantum/fig32.eps}

Fig 4.2.
\end{figure}

From now on we won't specify arrows any more, making the convention that 
they always point downwards at trivalent vertices.

Let us now start defining the basic data. The label set is
${\cal L}=\{m|$ $1\leq m<r\}$, the involution 
on ${\cal L}$ is the identity map, and the distinguished element is 1.

\medskip

If $D$ is a ce-disk, define $V(D,m)=V_m$ where
$V_m$ is the space of module homomorphisms ${\bf C}\rightarrow \ul{m}$.
By Schur's lemma $V_m=0$ unless $m=1$, and $V_1\tilde{=}{\bf C}$.
Let $\beta _1:{\bf C}\rightarrow \ul{1}$ be the isomorphism
$\lambda \rightarrow \lambda e_0$.

\medskip

If $A$ is a ce-annulus, $V(A,(m,n))=V^m_n=\{\phi :\ul{m}\rightarrow
\ul{n}$ $\phi$ module homomorphism$\}$. Like before, $V(A,(m,n))\neq 0$
if and only if $\ul{m}=\ul{n}$. In the latter case $V(A,(m,n))\tilde{=}
{\bf C}$. Let $\beta _m^m$ be the identity operator on $\ul{m}$.

\medskip

If $P$ is a ce-pair of pants define $V(P,(p,m,n))=V_p^{mn}$.
It follows that $V(P,(p,m,n))\neq 0$ if and only 
if $|m-n|+1\leq p \leq min\{m+n-1,2r-2-m-n\}$, in 
which case the vector space is isomorphic
to ${\bf C}$. 
An observation that is very useful in computations is the fact that
$V^{mn}_p\neq 0$ implies that $m+n+p$ is an odd integer, and $m-1+n-1+p-1$
is even.
Denote by $\beta ^{mn}_p$ one of the two projections in $V_p^{mn}$.
 One can choose these
projections such that $R\beta ^{mn}_{p}=\beta ^{pm}_n$
 where $R$ 
will be defined below. Choose also $\beta _1^{mm}=\psi _m$. Then 
$\beta _m^{1m}=R\beta _1^{mm}$ has the 
property that $\beta _m^{1m}(x\otimes e_0)=x$, for any $x\in \ul{m}$.

\medskip

The dual of $V_p^{mn}$ can be identified with $V_p^{nm}$.
Let $\psi :V_p^{mn}\rightarrow V_p^{nm}$ be given by
the diagram in Fig. 4.3.

\begin{figure}[htbp]
\centering
\leavevmode
\epsfxsize=3.4in
\epsfysize=1in
\epsfbox{Quantum/fig33.eps}

Fig. 4.3.
\end{figure}

The pairing  $<$ , $>:V_p^{mn}\bigotimes V_p^{nm}\rightarrow
{\bf C}$ is described in Fig. 4.4.

\begin{figure}[htbp]
\centering
\leavevmode
\epsfxsize=4.9in
\epsfysize=1.1in
\epsfbox{Quantum/fig34.eps}

Fig. 4.4.
\end{figure}

\ul{Remarks} 1. $<\alpha ,\beta >=X^2/(\sqrt{[m]}\sqrt{[n]}
\sqrt{[p]})<\alpha ,\psi (\beta )>_t$.

2. The pairing that we have defined is not quite the canonical
pairing between a space and its dual, it is rather the canonical pairing 
multiplied by a constant. We could have included that constant in the
definition of $\psi $, the way it is done in [FK], but we won't do it 
since it would complicate some computations.

\medskip

\ul{Proposition 4.1.} If $\alpha \in V_p^{mn}$ and $\beta \in V_p^{nm}$ then
$<\alpha ,\beta >=<\beta ,\alpha >$.

\medskip

\ul{Proof:} Since the space $V_p^{nm}$ 
is one dimensional we can assume that $\beta
=c\cdot \alpha \circ {\mathaccent20{R}}$, where $c$ is some complex number.
Using the fact that the operators $D$ can be pulled through crossings, we get
 $<\alpha ,\beta >=c<\alpha ,\alpha \circ {\mathaccent20{R}}>=$

\begin{figure}[htbp]
\centering
\leavevmode
\epsfxsize=6.1in
\epsfysize=1.3in
\epsfbox{Quantum/fig35.eps}

Fig. 4.5.
\end{figure}

$=c<\alpha \circ {\mathaccent20{R}}, \alpha>$. 
$\Box$

The dual of $V^m_m$ is $V_m^m$, $\psi =id$ and the pairing is given by
$<\beta _m^m,\beta _m^m>=X/[m]$, the dual of $V_1$ is also $V_1$,
$\psi $ is the identity and the pairing is given by $<\beta _1,
\beta _1>=1$.

In what follows, for an lce-morphism $f$, we will denote $V(f)$
also by $f$ whenever this seems unlikely to cause confusion.

 We now define the functor for the elementary moves that transform one
$DB$ structure into another. Exactly like in the case of vector spaces,
instead of describing the isomorphisms between vector spaces one can
describe the changes of basis, and then all isomorphisms can be regarded
as having the matrix equal to identity, with apropriate basis choice
for domain and range.
 Also in order to get positive crossings for the diagrams of morphisms
we  have to consider negative crossings in the diagrams below (see the
usual connection between changes of basis and isomorphisms). 

If $\sigma $ is an lce-surface and $K$ is the move defined in the 
first section applied to a band along a circle labeled by $m$,
define
$V(K):V(\sigma )\rightarrow V(\sigma )$
to be the multiplication by $(-1)^{m-1}$.

\medskip

If $P$ is an lce-pair of pants, with boundary components labeled 
respectively $p,m,n$ define $T_1$, $B_{23}$ and $R$ like in Fig. 4.6.
 
\begin{figure}[htbp]
\centering
\leavevmode
\epsfxsize=4.2in
\epsfysize=2.5in
\epsfbox{Quantum/fig36.eps}

Fig 4.6.
\end{figure}

For two pairs of pants glued along a boundary component
we let 
\begin{eqnarray*}
F:\bigoplus_{p\in {\cal L}}
V_p^{mn}\otimes V_p^{kl}\rightarrow \bigotimes_{q\in {\cal L}}V_p^{lm}
\otimes V_p^{nk}
\end{eqnarray*}
be defined by the  relation in Fig. 4.7.

\begin{figure}[htbp]
\centering
\leavevmode
\epsfxsize=5in
\epsfysize=1.5in
\epsfbox{Quantum/fig37.eps}

Fig. 4.7.
\end{figure}

 For an lce-pair of pants, with the bottom boundary components
glued together we define $S:\bigoplus_{m\in {\cal L}}V_p^{mm}\rightarrow
\bigoplus_{n\in {\cal L}}V_p^{nn}$, by the diagram in Fig. 4.8.

\begin{figure}[htbp]
\centering
\leavevmode
\epsfxsize=3.6in
\epsfysize=1.1in
\epsfbox{Quantum/fig38.eps}

Fig. 4.8.
\end{figure}

Also
$P_{12}:V_p^{mn}\bigotimes V_p^{kl}\rightarrow V_p^{kl}\bigotimes V_p^{mn},
\ \ P(\alpha \otimes \beta )=(\beta \otimes \alpha )$\\
\ \ \ \ \ \ $A:V_p^{mn}\bigotimes V_m^{m}\rightarrow V_p^{mn},
\ \ A(\beta ^{mn}_p\otimes
\beta ^m_m)=\beta ^{mn}_p$\\
\ \ \ \ \ \ $D:V_m^{1m}\bigotimes V_1\rightarrow 
V_m^m,\ \ D(\beta ^{1m}_m\otimes
\beta _1)=\beta _m^m$\\
\ \ \ \ \ \ $A:V_m^m\bigotimes V_m^m\rightarrow V_m^m, 
\ \ A( \beta ^{m}_m\otimes
\beta _m^m)=\beta _m^m$\\
\ \ \ \ \ \ $D:V_1^1\bigotimes V_1\rightarrow V_1, \ \ D(\beta ^{1}_1\otimes
\beta _1)=\beta _1$

The map $C=V((id,1)):V(\sigma )\rightarrow V(\sigma )$ is the multiplication by
$exp(3\pi (r-2)i/(4r))$, and finally $S:{\cal L}\rightarrow {\bf C}$
is given by $S(m)=[m]/X$.

\medskip

\ul{Remark} These maps are canonical, and they can be found in [Wa] and
[FK]. What is important in our situation is the exact location of the 
coupons. Let us also note that the elements of the matrix of $F$ are slight
modifications of the $6j$-symbols.

\bigskip
\begin{center}
{\bf 5. RELATIONS THAT THE BASIC DATA MUST SATISFY}
\end{center}
\bigskip
  
The first result exhibits the Moore-Seiberg equations that
a modular functor on the category of ce-surfaces must satisfy. It
is the analogue for ce-surfaces  of  Theorem 6.4 from [Wa].

\bigskip

\ul{Theorem 5.1.} A basic data determines a modular functor $V$ satisfying the 
axioms (3.1)--(3.10) (a modular functor) if and only if the following
relations hold:

1. relations at the level of a ce-pair of pants:

\ a) $T_1B_{23}=B_{23}T_1$,\ $T_2B_{23}=B_{23}T_3$, \ $T_3B_{23}=B_{23}T_2$,
where $T_2=RT_1R^{-1}$ and $T_3=R^{-1}T_1R$.

\ b) $B_{23}^2=T_1T_2^{-1}T_3^{-1}$

\ c) $R^3=1$

\ d) $RB_{23}R^2B_{23}RB_{23}R^2=B_{23}RB_{23}R^2B_{23}$

\ e) $K^2=1$

2. relations defining the inverses of $F$ and $S$:

\ a) $P_{12}K_1^{(1)}F^2=1$

\ b) $K_2T_3^{-1}B_{23}^{-1}S^2=1$

3. relations coming from ``codimension 2 singularities'':

\ a) $K_1^{(1)}P_{13}R^{(2)}F^{(12)}K_1^{(1)}R^{(2)}K_1^{(2)}F^{(23)}R^{(2)}
F^{(12)}K_1^{(1)}R^{(2)}K_1^{(2)}F^{(23)}R^{(2)}F^{(12)}=1$

\ b) $T_3^{(1)}FB_{23}^{(1)}FB_{23}^{(1)}FB_{23}^{(1)}=1$

\ c) $C^{-1}K_2B_{23}^{-1}T_3^{-2}ST_3^{-1}ST_3^{-1}S=1$

\ d) $R^{(1)}(R^{(2)})^{-1}FS^{(1)}FB_{23}^{(2)}
B_{23}^{(1)}= K_1^{(1)}FS^{(2)} T_3^{(2)}(T_1^{(2)})^{-1}B_{23}^{(2)}
F$

4. relations involving annuli and disks:

\ a) $F(\beta _p^{mn}\otimes \beta _p^{p1})=\beta _m^{1m}\otimes \beta _m^{np}$

\ b) $A^{(12)}_2D_3^{(13)}=
D_2D_3^{(13)}$

\ c) $A^{(12)}A^{(23)}=
A^{(23)}A^{(12)}$

5. relations coming from duality:

\ -- for any elementary move $f$, one has $f^+=\bar{f}$ where $f^+$ is
the dual of $f$ with respect to the $<$ , $>$ pairing, and $\bar{f}$ is
the morphism induced by $f$ on $-\sigma $,

6. relations expressing the compatibility between the pairing, and moves 
$A$ and $D$:

\ a) $<\beta _m^m,\beta _m^m>=S(m)^{-1}$

\ b) $<\beta _m^{m1},\beta _m^{m1}>=S(1)^{-1}S(m)^{-1}$.

\bigskip

The theorem needs some clarifications.
The first group of relations holds on a lce-pair of pants. For the rest of 
the relations the superscripts indicate the elementary
surfaces to which the move is applied, while the  subscripts
indicate the number of the boundary component. The surfaces on which these
relations hold are the ones in Fig. 5.1,
(where we allow the bands along the boundary components
to be indexed by any elements of the Klein group).

\begin{figure}[htbp]
\centering
\leavevmode
\epsfxsize=4.6in
\epsfysize=5.8in
\epsfbox{Quantum/fig39.eps}

Fig. 5.1.
\end{figure}

\bigskip

\ul{Remark} As the reader can see, there are no listed relations that
are to be satisfied at the level of a torus. This is because
those relations are redundant, the torus can be obtained by capping the
punctured torus with a disk.

\medskip

\ul{Proof of the theorem:} It is not hard to see that the theorem is a
consequence of Theorem 6.4 in [Wa]. For a better understanding let us
start by briefly explaining the origin of the Moore-Seiberg  relations.
 
The first group of relations are satisfied at the level of a pair of pants.
If one forgets about the bands, one gets the presentation of the mapping class
group of a pair of pants. Since the mapping class group is isomorphic
to the group of moves (namely of transformations of seams and numbers),
 the same 
relations give a presentation of the latter group. Lifting these relations 
to ce-pairs of pants one gets 1. a)--e).

The second and the third group of relations can be obtained via Cerf
theory [Ce] by using the techniques described in [HT]. The first group
comes from cancellations of consecutive crossings in the graph of the
height function, and give the inverses of $F$ and $S$. The second
group arises from triangle singularities. Actually, the triangle singularities
give rise to a larger number of relations, but as remarked in [Wa],
they are all consequences of four fundamental ones. These four 
relations are liftings at the level of ce-surfaces
of the following relations:
the pentagon, shown in Fig.  5.2, the F-triangle, the S-triangle, both 
shown in Fig. 5.3 and the (FSF)$^2$-cell shown in Fig. 5.4.
  
\begin{figure}[htbp]
\centering
\leavevmode
\epsfxsize=4.5in
\epsfysize=2.2in
\epsfbox{Quantum/fig40.eps}

Fig. 5.2.
\end{figure}

The fourth group of relations expresses the fact that the operation of 
cancelling an annulus or a disk is preserved by the functor. Finally,
the last two groups of relations must hold in order for the duality
axiom to be satisfied.

Let us now proceed to proving that these relations are sufficient in
order for the basic data to define a modular functor.

\begin{figure}[htbp]
\centering
\leavevmode
\epsfxsize=5.7in
\epsfysize=2.3in
\epsfbox{Quantum/fig41.eps}

Fig. 5.3.
\end{figure}

Given a ce-surface $\sigma$ let us consider the CW complex $\Lambda$,
whose vertices are all the ce-surfaces that can be obtained from
$\sigma $ by performing some moves; whose edges are elementary moves
between these surfaces, and whose 2-cells are the groups of relations
1, 2, 3, 4, together with all the 2-cells which express
disjoint commutativity between elementary moves (which are obviously
satisfied by any basic data). The fact that $\Lambda $ is 
connected follows from the construction. Let us prove that it is
simply connected. 

If we consider the CW complex $\Gamma$ obtained from $\Lambda$
by forgetting the bands (i.e. the analogous CW complex for
DAP-decompositions), then Theorem 6.4 in [Wa] asserts that 
$\Gamma $ is simply connected. 

Let $\phi :\Lambda \rightarrow \Gamma$
be the quotient function. We see that for any 0-cell $\sigma _0$
of $\Gamma$, $\phi ^{-1}(\phi (\sigma _0))$ consists out of the  subcomplex
of $\Lambda $ whose 0-cells are the ce-surfaces that can be 
obtained from $\sigma _0$ by performing moves of type $K$. 
Since any two such moves commute, $\phi ^{-1}(\phi (\sigma _0))$ is
simply connected. It follows that $\Lambda$ is also simply connected, 
hence the relations exhibited in the statement are sufficient for the 
functor $V$ to be well defined. 

\begin{figure}
\centering
\leavevmode
\epsfxsize=5in
\epsfysize=2in
\epsfbox{Quantum/fig42.eps}

Fig. 5.4.
\end{figure}

The rest of the relations imply that the   
duality axiom for $V$ holds, and the  theorem is proved.$\Box$

\medskip

The next result gives necessary and sufficient conditions for the
partition function $Z$ to be well defined.

\bigskip

\ul{Theorem 5.2.} Given  a functor
$V$ that satisfies the relations from Theorem 1 the partition 
function that it determines is well defined if and only if the following
two conditions hold:

\medskip

\hspace{40mm} $a)$ $S(m)=S_{1m}$ 
where $[S_{xy}]_{x,y}$  is the matrix of move $S$,
\begin{eqnarray*}
\ b) F(\beta _1^{mm}\otimes\beta _1^{nn})=
\sum_p\frac{(-1)^{n-1}id_{mnp}}{S(m)S(n)}
\end{eqnarray*}
where $id_{mnp}$ is the identity matrix in $(V_p^{mn})^*\otimes V_p^{mn}$.

\medskip

\ul{Proof:} Let us first convince ourselves that the two relations are
necessary. 

a) Let us view a ball as the mapping cylinder of the identity
map of the disk. Then the invariant of the ce-ball given in Fig. 5.5.

\begin{figure}[htbp]
\centering
\leavevmode
\epsfxsize=1.2in
\epsfysize=1.2in
\epsfbox{Quantum/fig43.eps}
2~
Fig. 5.5.
\end{figure}

\noindent is $\beta _1\otimes \beta _1$, by axiom (3.10).
Let us consider the  chain of transformations from Fig. 5.6,
where we write under each extended ce-3-manifold its invariant.

\begin{figure}[htbp]
\centering
\leavevmode
\epsfxsize=4.5in
\epsfysize=2.3in
\epsfbox{Quantum/fig44.eps}

Fig. 5.6.
\end{figure}

\noindent where $S_T$ is the composition of elementary moves 
described in Fig. 1.16 ( the $S$-move on a torus). By the identity described in
Fig. 1.18 the latter ce-manifold is equal to the one from 

\begin{figure}[htbp]
\centering
\leavevmode
\epsfxsize=2.3in
\epsfysize=1.1in
\epsfbox{Quantum/fig45.eps}

Fig. 5.7.
\end{figure}

\noindent Fig. 5.7. which is the mapping cylinder of the identity map on a
ce-annulus. The invariant of this ce-manifold is $\bigoplus
_m(id:V_m^m\rightarrow V_m^m)$, and from Theorem 5.1 relation 6. a)
it follows that this is equal 
to $\bigoplus _mS(m)\beta _m^m\otimes \beta _m^m$.
This gives $S(m)=S_{1m}$, $\forall m$.

b) Similarly we have the following chain of transformations
from Fig. 5.8 followed by that from Fig. 5.9. 

\begin{figure}[htbp]
\centering
\leavevmode
\epsfxsize=5.3in
\epsfysize=0.9in
\epsfbox{Quantum/fig46.eps}

\centering
\leavevmode
\epsfxsize=5.3in
\epsfysize=0.9in
\epsfbox{Quantum/fig47.eps}

Fig. 5.8.
\end{figure}

To it corresponds the following chain of transformations of the invariants:

\hspace{10mm}$\beta _1^{1}\otimes \beta _1^1 
\rightarrow \beta _1^{11}\otimes 
\beta _1\otimes \beta _1^{11}\otimes \beta _1 
\rightarrow \beta _1^{11}\otimes \beta _1\otimes \beta _1^{11}
\otimes \beta _1\rightarrow \\
\bigoplus _{m,n}S(m)S(n)\beta _1^{mm}\otimes \beta _1\otimes 
\beta _1^{nn}\otimes \beta _1
\rightarrow \bigoplus _{m,n}(-1)^{2m}S(m)S(n)\beta _1^{mm}
\otimes \beta _1\otimes 
\beta _1^{nn}\otimes \beta _1 \rightarrow \\
\bigoplus _{m,n}(-1)^{2m}S(m)S(n)\beta _1^{mm}\otimes 
\beta _1^{nn}
\rightarrow \bigoplus _{m,n}(-1)^{2n}S(m)S(n)F(\beta _1^{mm}\otimes 
\beta _1^{nn})$

\begin{figure}[htbp]
\centering
\leavevmode
\epsfxsize=5.3in
\epsfysize=1in
\epsfbox{Quantum/fig48.eps}

\centering
\leavevmode
\epsfxsize=3.8in
\epsfysize=1.2in
\epsfbox{Quantum/fig49.eps}

Fig. 5.9.
\end{figure}

\noindent This is the mapping cylinder of a ce-pair of pants, so
relation b) is implied by the mapping cylinder axiom.

For the proof of sufficiency, let $\mu =(M,D,B,n)$ be a ce-3-manifold. 
The manifold $M$ can be obtained by successively attaching 0, 1, 2 or
3-handles to a ball. By putting DB-structures on the handles, (and eventually
performing some moves), we  see that one can obtain in this way a 
ce-3-manifold $(M,D',B',n')$. The right choice of framing for the ce-handles
gives $n=n'$. Since $\partial M$  is a closed surface, there is a move 
transforming $D'$ and $B'$ into $D$ and $B$. Thus $\mu $ can be obtained 
by successively adding handles to a ce-ball, and eventually performing moves
on the boundary.

As a consequence we get that the techniques of Morse theory described in [Wa] 
work mutatis mutandis to show that the two relations in the statement imply
that $Z$ is well defined.$\Box$
 
\bigskip
\begin{center}
{\bf 6. THE VERIFICATION OF THE COMPATIBILITY CONDITIONS FOR
THE BASIC DATA}
\end{center}
\bigskip
  
In this section we will check that our basic data satisfies the conditions
from Theorems 5.1 and 5.2.
We start with the Moore-Seiberg relations. Since most of the proofs are 
similar to those from [FK] we will only check one relation from
each group, for the rest we refer the reader to [FK]. 
The first relation we want to check is 1. c).

\medskip

\ul{Lemma 6.1.} If $\alpha \in V_p^{mn}$ and $\beta \in V_p^{nm}$
then $<R\alpha ,\beta >=(-1)^{m-1}<\alpha ,R\beta >$.

\medskip

\ul{Proof:} The proof is described in Fig. 6.1.$\Box$

\begin{figure}[htbp]
\centering
\leavevmode
\epsfxsize=6in
\epsfysize=3.9in
\epsfbox{Quantum/fig50.eps}
 
Fig. 6.1.
\end{figure}

Applying the lemma we get that for $\alpha \in V_p^{mn}$ and 
$\beta \in V_p^{nm}$, \ $<\alpha , R^3\beta >=
(-1)^{m-1}<R\alpha ,R^2\beta >$. From here the proof proceeds like 
in Fig. 6.2.
where the last equality follows from Proposition 4.3.$\Box$

\begin{figure}
\centering
\leavevmode
\epsfxsize=5.4in
\epsfysize=2.8in
\epsfbox{Quantum/fig51.eps}

Fig. 6.2.
\end{figure}

\medskip

\ul{Corollary 6.2.} The  identities from Fig. 6.3 hold.

\begin{figure}[htbp]
\centering
\leavevmode
\epsfxsize=4.5in
\epsfysize=2.5in
\epsfbox{Quantum/fig52.eps}

Fig.6.3.
\end{figure}

\ul{Proof:} The first identity follows by applying $\psi ^{-1}$
to both terms of  the equality $R^3\alpha =\alpha $. The second
identity follows from the first by taking the adjoint.$\Box$

\bigskip

\ul{Proposition 6.3.} Assume that $C_i:\ul{1}
\rightarrow \ul{m_i}\otimes \ul{n_i}
\otimes \ul{m_{i+1}}\otimes \ul{n_{i+1}}$, $i=1,\cdots ,N$, $N+1=1$
are morphisms of representations. Let $\beta _q^i=\beta _q^{m_in_i}$,
$q\in {\cal L}$. Then the identity  from Fig. 6.4 holds.

\begin{figure}[htbp]
\centering
\leavevmode
\epsfxsize=4.5in
\epsfysize=1.9in 

\epsfbox{Quantum/fig53.eps}

Fig. 6.4.
\end{figure}

\ul{Proof:} We proceed by transforming the left hand side like in
Fig. 6.5. 

\begin{figure}[htbp]
\centering
\leavevmode
\epsfxsize=4.8in
\epsfysize=2.8in
\epsfbox{Quantum/fig54.eps}

Fig. 6.5.
\end{figure}

Since by Schur's lemma all the factors in the middle are multiples 
of the identity on $q$ they can be inserted into the first factor to get
the expression from Fig. 6.6, which is further equal to the one in 
Fig. 6.7 since in the first factor of this latter expression
the only thing that matters is the copy of $q$ in each $\ul{m_i}\otimes
\ul{n_i}$.

\begin{figure}[htbp]
\centering
\leavevmode
\epsfxsize=5.3in
\epsfysize=1.1in
\epsfbox{Quantum/fig55.eps}

Fig. 6.6
\end{figure}

\begin{figure}[htbp]
\centering
\leavevmode
\epsfxsize=5.2in
\epsfysize=1.1in
\epsfbox{Quantum/fig56.eps}

Fig. 6.7.

\end{figure}

Further, this is equal to 
the diagram from Fig. 6.8
which by Proposition 7 in [FK] is equal to the desired expression.
$\Box$

\begin{figure}[htbp]
\centering
\leavevmode
\epsfxsize=5.2in
\epsfysize=1.8in
\epsfbox{Quantum/fig57.eps}

Fig. 6.8.
\end{figure}

\medskip

We are now ready to proceed with the proof of relation 2. a). Write it as 
$FP_{12}K_1^{(1)}F=1$. Let $\alpha \otimes \beta \in V_p^{ij}\bigotimes
V_p^{kl}$ and $\delta \otimes \gamma \in V_q^{ij}\bigotimes V_q^{kl}$.
Recall that for $u,v,w\in {\cal L}$, $\beta _w^{uv}\in V_w^{uv}$
is the orthogonal projection. We have

\medskip

\ \ $<FP_{12}K_1^{(1)}F(\alpha \otimes \beta ),\delta \otimes \gamma >_t=
<F( \alpha \otimes \beta ), (FP_{12}K_1^{(1)})^*(\delta \otimes \gamma )>_t=$

\medskip

$\sum_r <F\alpha \otimes \beta , \beta _r^{li}\otimes \beta _r^{jk}>_t
<\beta _r^{li}\otimes \beta _r ^{jk}, (FP_{12}K_1^{(1)})^*(\delta
\otimes \gamma )>_t=$

\medskip

$\sum_r <F\alpha \otimes \beta , \beta _r^{li}\otimes \beta _r^{jk}>_t
<FK_1^{(1)}\beta _r^{jk}\otimes \beta _r ^{li}, \delta
\otimes \gamma >_t$. 

\medskip

From here we proceed like in Fig. 6.9.
Using Proposition 6.3, we see that we can continue our computation like
in
Fig. 6.10.

\begin{figure}[htbp]
\centering
\leavevmode
\epsfxsize=5.2in
\epsfysize=2.3in
\epsfbox{Quantum/fig58.eps}

Fig. 6.9.
\end{figure}

\begin{figure}[htbp]
\centering
\leavevmode
\epsfxsize=4.6in
\epsfysize=3.3in
\epsfbox{Quantum/fig59.eps}

Fig. 6.10.
\end{figure}

By Schur's lemma and the definition of the trace
pairing this is equal to
\begin{eqnarray*}
\frac{1}{\sqrt{[p]}\sqrt{[q]}}\delta _{p,q}<\alpha ,\delta>_t
<\beta ,\gamma>_t [p]
\end{eqnarray*}
 where $\delta _{p,q}$ is the Kronecker symbol. The latter
is equal to $\delta _{p,q}<\alpha
,\delta>_t
<\beta ,\gamma>_t$ which shows that the matrix of our morphism
is the identity matrix.$\Box$

\medskip

Let us now prove the pentagon identity. Rewrite it by using
1. c) and 2. a) in the form
\begin{eqnarray*}
R^{(2)}F^{(12)}K_1^{(1)}R^{(2)}K_1^{(2)}F{^(23)}R^{(2)}F^{(12)}=
F^{(23)}P_{23}P_{12}(R^{(1)})^2F^{(12)}(R^{(2)})^2P_{13}K_1^{(1)}.
\end{eqnarray*}
We prove the identity by checking the action on projectors and by
using the trace pairing. For the left hand side we have

\ $<R^{(2)}F^{(12)}K_1^{(1)}R^{(2)}K_1^{(2)}F^{(23)}R^{(2)}F^{(12)}
\beta _p^{ij}
\otimes \beta _p^{qk}\otimes \beta _q^{lm},\beta _r^{mk}\otimes \beta _s^{ri}
\otimes \beta _s^{jl}>_t=$

\ $<F^{(12)}K_1^{(1)}R^{(2)}K_1^{(2)}F^{(23)}R^{(2)}F^{(12)}\beta _p^{ij}
\otimes \beta _p^{qk}\otimes \beta _q^{lm},\beta _r^{mk}\otimes \beta _r^{is}
\otimes \beta _s^{jl}>_t=$

 \medskip

\ $\sum_t<F\beta _p^{ij}\otimes \beta _p^{qk}, \beta _t^{ki}\otimes \beta
_t^{jq}>_t<FR^{(1)}\beta _t^{jq}\otimes \beta _q^{lm},\beta _s^{mt}
\otimes \beta _s^{jl}>_t\cdot$

$<FK_1^{(1)}R^{(2)}K_1^{(2)}\beta _t^{ki}\otimes \beta _s^{mt},
\beta _r^{mk}\otimes \beta _r^{is}>_t=$

\medskip

\ $\sum_t(-1)^{2m}<F\beta _p^{ij}\otimes \beta _p^{qk}, 
\beta _t^{ki}\otimes \beta
_t^{jq}>_t<FR^{(1)}\beta _q^{tj}\otimes \beta _q^{lm},\beta _s^{mt}
\otimes \beta _s^{jl}>_t\cdot$

$<FK_1^{(1)}R^{(2)}K_1^{(2)}\beta _t^{ki}\otimes \beta _t^{sm},
\beta _r^{mk}\otimes \beta _r^{is}>_t$

From here the computation continues like in Fig 6.11.

\begin{figure}[htbp]
\centering
\leavevmode
\epsfxsize=5.5in
\epsfysize=2.7in
\epsfbox{Quantum/fig60.eps}

Fig. 6.11.
\end{figure}

Further we do two rotations in the second factor and a flip in the third to get
the expression from Fig. 6.12.
By using 1. b) in the third factor we get the diagram from Fig. 6.13.
Using  Proposition 6.3 we see that we can continue like in Fig. 6.14.

\begin{figure}
\centering
\leavevmode
\epsfxsize=5.5in
\epsfysize=2.7in
\epsfbox{Quantum/fig61.eps}

Fig. 6.12.
\end{figure}

\begin{figure}
\centering
\leavevmode
\epsfxsize=5.5in
\epsfysize=2.8in
\epsfbox{Quantum/fig62.eps}

Fig. 6.13.
\end{figure}

\begin{figure}[htbp]
\centering
\leavevmode
\epsfxsize=5.9in
\epsfysize=5.4in
\epsfbox{Quantum/fig63.eps}

Fig. 6.14.
\end{figure}

For the right hand side we have

\medskip

$<F^{(23)}P_{23}P_{12}(R^{(1)})^2F^{(12)}
(R^{(2)})^2P_{13}K_1^{(1)}\beta _p^{ij}
\otimes \beta _p^{qk}\otimes \beta _q^{lm},\beta _r^{mk}\otimes \beta _s^{ri}
\otimes \beta _s^{jl}>_t=$

\medskip

$(-1)^{2p}<F\beta _q^{lm}\otimes \beta _q^{kp}, \beta _r^{pl}\otimes 
\beta _r^{mk}>_t<F\beta _p^{ij}\otimes \beta _p^{lr}, \beta _s^{ri}\otimes
\beta _s^{jl}>_t.$

This is equal  to the morphism described by the diagram in Fig 6.15.

\begin{figure}[htbp]
\centering
\leavevmode
\epsfxsize=5.8in
\epsfysize=1.2in
\epsfbox{Quantum/fig64.eps}

Fig. 6.15.
\end{figure}

By applying one rotation in the first factor and other two rotations
in the second, we get 
the expression from Fig. 6.16.

\begin{figure}[htbp]
\centering
\leavevmode
\epsfxsize=5.3in
\epsfysize=2.4in
\epsfbox{Quantum/fig65.eps}

Fig. 6.16.
\end{figure}

We do  a flip in the second factor use 1. c)
and continue like in Fig. 6.17.

\begin{figure}[hpbt]
\centering
\leavevmode
\epsfxsize=5.5in
\epsfysize=2.7in
\epsfbox{Quantum/fig66.eps}

\centering
\leavevmode
\epsfxsize=5.5in
\epsfysize=2.5in
\epsfbox{Quantum/fig67.eps}

Fig. 6.17.
\end{figure}

From here, after   using Proposition 7. c) 
in [FK] we continue wth the computations like in Fig. 6.18 to finally get
the expression from Fig. 6.19, and the identity is proved.

\begin{figure}[hpbt]

\centering
\leavevmode
\epsfxsize=5.9in
\epsfysize=4.3in
\epsfbox{Quantum/fig68.eps}

Fig. 6.18.
\end{figure}

\begin{figure}[htbp]
\centering
\leavevmode
\epsfxsize=5.2in
\epsfysize=2in
\epsfbox{Quantum/fig69.eps}

Fig. 6.19.
\end{figure}

For the relation 4. a) let  us first remark that $F(\beta _p^{mn}\otimes
\beta _p^{p1})\in V_m^{1m}\bigotimes V_m^{np}$, so we only have to 
compute $<F(\beta _p^{mn}\otimes
\beta _p^{p1}), \beta _m^{1m}\otimes \beta _m^{np}>_t$. This is done in
Fig. 6.20. 

\begin{figure}[htbp]
\centering
\leavevmode
\epsfxsize=5.4in
\epsfysize=3.7in
\epsfbox{Quantum/fig70.eps}

Fig. 6.20.
\end{figure}

\medskip

The fifth group of identities is satisfied because of the following equalities

\ \ $\bar{B}_{23}=B_{23}^{-1}$, $\bar{K}_i=K_i$, $\bar{T}_i=T_i^{-1}$,
$\bar{R}=R^{-1}K_1K_3$, $\bar{S}=S^{-1}$ and $\bar{F}=F^{-1}$.

\medskip

The last group of relations is clearly satisfied.$\Box$

Let us proceed in verifying that the two relations from Theorem 5.2 hold.
The first one holds by definition. For the second one
let us first remark that the identity matrix on $V^{mn}_p$ has the
form

\newpage

 $\bigoplus _{m,n,p}\hat{\beta _{p}^{mn}}\otimes \beta _p^{mn}$,
where  $\hat{\beta _{p}^{mn}}\in V^{nm}_p$ is 
the base elements dual to $\beta _p^{mn}$
with respect to the pairing. 
We only have to check that 
\begin{eqnarray*}
<F\beta ^{mm}_1\otimes \beta ^{nn}_1,\beta ^{mn}_p\otimes \beta ^{nm}_p>=
\frac{(-1)^{n-1}<\beta ^{mn}_p,\beta ^{nm}_p>}{S(m)S(n)}.
\end{eqnarray*}

We have
\begin{eqnarray*}
<F\beta ^{mm}_1\otimes \beta ^{nn}_1,\beta ^{mn}_p\otimes \beta ^{nm}_p>=
\frac{X^4}{[m][n][p]}<F\beta ^{mm}_1
\otimes \beta ^{nn}_1,\psi (\beta ^{mn}_p)\otimes \psi (\beta ^{nm}_p)>_t
\end{eqnarray*}

The value of this expression is given in Fig. 6.21.

\begin{figure}[htbp]
\centering
\leavevmode
\epsfxsize=4in
\epsfysize=1.2in
\epsfbox{Quantum/fig71.eps}

Fig. 6.21.
\end{figure}

By using the identity from Fig. 6.22, we see that we can continue like in Fig.
6.23 to get the desired result.$\Box$

\begin{figure}[htbp]
\centering
\leavevmode
\epsfxsize=3in
\epsfysize=0.6in
\epsfbox{Quantum/fig72.eps}

Fig. 6.22.
\end{figure}

\begin{figure}[htbp]

\centering
\leavevmode
\epsfxsize=5in
\epsfysize=2.4in
\epsfbox{Quantum/fig73.eps}

Fig. 6.23.
\end{figure}

\bigskip

\newpage

{\bf REFERENCES}

\bigskip

[A] Atiyah, M. F., {\em The Geometry and Physics of knots},
Lezioni Lincee, Accademia Nationale de Lincei, Cambridge Univ.
Press, 1990.

[Ce] Cerf, J., { \em La stratification naturelle et le th{\'{e}}oreme de la
pseudo-isotopie}, Publ. Math. I.H.E.S., {\bf 39}(1969), 5--173.

[FK] Frohman, Ch., Kania-Bartoszynska, J., {\em $SO(3)$ topological quantum
field theory}, preprint, 1994.

[G] Gelca, R., {\em The quantum invariant of the complement of a
link}, preprint.

[HT] Hatcher, A., Thurston, W., {\em A presentation of the mapping 
class group of a closed orientable surface}, Topology, {\bf 19}(1980),
221--237.

[J] Jones, V., F., R., {\em Polynomial invariants of knots via von
Neumann algebras}, Bull. Amer. Math. Soc., {\bf 12}(1998), 103--111.

[KM] Kirby, R., Melvin, P., {\em The 3-manifold invariants of
Witten and Reshetikhin--Turaev for $sl(2,{\bf C})$}, Inventiones Math.,
{\bf 105}(1991), 547--597.

[RT] Reshetikhin, N. Yu., Turaev, V. G., {\em Invariants of 3-manifolds
via link polynomials and quantum groups}, Inventiones Math., {\bf 103}(1991),
547--597.

[T] Turaev, V., G., {\em Quantum invariants of Knots and 3-manifolds},
de Gruyter Studies in Mathematics, de Gruyter, Berlin--New York, 1994. 

[Wa] Walker, K., {\em On Witten's 3-manifold invariants}, preprint, 1991. 

[W] Wall, C. T. C., {\em Non-additivity of the signature}, Inventiones Math.,
{\bf 7}(1969), 269--274.

\bigskip

Department of Mathematics, The University of Iowa, Iowa City, IA
52242 (mailing address)
{\em E-mail: rgelca@math.uiowa.edu}
\medskip

and

\medskip

Institute of Mathematics of Romanian Academy, P.O.Box 1-764,
70700 Bucharest, Romania.

\end{document}